 \documentclass[12pt, twoside]{amsart}

\usepackage{amsmath,amssymb,amsthm,amscd, amsxtra, amsfonts}

\makeatletter \@addtoreset{equation}{section}

\makeatother \makeatletter

\def \<{\langle}
\def \>{\rangle}

\def \a{\alpha }

\def \l{\lambda }

\newtheorem{theorem}{Theorem}[section]

\newtheorem{lemma}{Lemma}[section]

\newtheorem{remark}{Remark}[section]
\newtheorem{definition}{Definition}[section]
\newtheorem{proposition}{Proposition}[section]

\newcommand{\bea}{\begin{eqnarray}}
\newcommand{\eea}{\end{eqnarray}}
\newcommand{\be}{\begin {equation}}
\newcommand{\ee}{\end{equation}}

\newcommand{\h}{\frak h}
\newcommand{\wt}{{\rm {wt} }   }

\newcommand{\Z}{\Bbb Z}
\newcommand{\W}{\mathcal W}

\newcommand{\Zp}{{\Bbb Z}_{>0} }

\newcommand{\N}{{\Bbb Z}_{\ge 0} }
\newcommand{\C}{\Bbb C}

\newcommand{\1}{\bf 1}

\newcommand{\la}{\langle}
\newcommand{\ra}{\rangle}

\begin{document}

\title[]{ Classification of irreducible modules of certain subalgebras of
  free boson vertex algebra}
 \footnote{1991 Mathematics Subject  Classification.
Primary 17B69; Secondary 17B68, 81R10.}
\author{ Dra\v zen Adamovi\' c }
\address{  Department of Mathematics, University of Zagreb  \\
Bijeni\v cka 30, 10 000 Zagreb, Croatia  \\ e-mail: \
adamovic@math.hr }
 \maketitle

\begin{abstract}
Let $M(1)$ be the  vertex  algebra for a  single free boson. We
classify irreducible modules of certain vertex subalgebras of
$M(1)$ generated by two generators. These subalgebras  correspond
to  the ${\W}(2, 2 p -1)$--algebras with central charge $ 1- 6
\frac{(p - 1) ^{2}}{p} $ where $p$ is a positive integer, $p \ge
2$. We also determine the associated Zhu's algebras.
\end{abstract}

\section{ Introduction}

Let $M(1)$ be the  vertex  algebra generated by a single free
boson. For every $z \in {\C}$, this vertex algebra contains a
Virasoro vertex operator subalgebra with central charge $1 - 12 z
^{2}$ (cf. \cite{MN}, \cite{K}). Therefore, $M(1)$ can be treated
as a vertex operator algebra of rank $1 - 12 z ^{2}$.

 The vertex operator algebra $M(1)$ has a family of
irreducible ${\N}$--graded (untwisted) modules $M(1, \l)$, $\l \in
{\C}$, and a $\Z_2$--twisted irreducible module $M(1) ^{\theta}$.
Some  subalgebras of $M(1)$ have  property that any irreducible
module for such subalgebra can be constructed from twisted or
untwisted modules for $M(1)$. In particular, this is true for the
orbifold vertex operator algebra $M(1) ^{+}$ (cf. \cite{DN}).
Another interesting example of such vertex operator algebra was
studied by W. Wang in \cite{W2}, \cite{W3}. He showed that
${\W}(2,3)$--algebra with central charge $-2$ can be realized as a
subalgebra of  $M(1)$,  and that every irreducible
${\W}(2,3)$--module is obtained  from $M(1)$--modules $M(1,\l)$.

Recall that for any vertex operator algebra $V$, Zhu in \cite{Z}
constructed an associative algebra $A(V)$ such that there is
one-to-one correspondence between the irreducible $V$--modules and
the irreducible $A(V)$--modules.
 In the cases mentioned above  the authors explicitly determine the corresponding
Zhu's algebras as certain quotients of the polynomial algebra
${\C}[x,y]$.

In the present paper we shall investigate certain vertex
subalgebras of $M(1)$ generated by two generators. Let us describe
these algebras. Let $\omega$ be the Virasoro element in $M(1)$
which generates a subalgebra isomorphic to the simple Virasoro
vertex operator algebra $W_0 \cong L(c_{p,1},0)$ (cf. \cite{W1},
\cite{FZ}) where $c_{ p,1} = 1 - 6 (p-1) ^{2} / p$, $p \ge 2$.
This vertex operator algebra can be extended by the primary vector
$H \in M(1)$ of level $2 p -1$ which is defined using Shur
polynomials. This new extended vertex operator algebra, which we
denote by $\overline{M(1)}$, is known in the physical literature
as ${\W}(2, 2 p-1)$--algebra (cf. \cite{EFH}, \cite{F}, \cite{H}).

We also study another description of the vertex operator algebra
$\overline{M(1)}$. It is  well-known (cf. \cite{FB}, \cite{FFr})
that some vertex subalgebras of $M(1)$ can be defined as kernels
of screening operators. In our particular case we consider two
screening operators $Q$ and $\widetilde{Q}$ defined using the
formalism of (generalized) vertex operator algebras (see Section
\ref{stwist}).  Then $\mbox{Ker}_{M(1)}Q$ is isomorphic to the
vertex operator algebra $L(c_{p,1},0)$, and $\mbox{Ker}_{M(1)}
\widetilde{Q}$ is isomorphic to $ \overline{M(1)}$. These
cohomological characterizations provide some deeper information on
the structure of $\overline{M(1)}$. It turns out that
$\overline{M(1)}$ is a completely reducible module for the
Virasoro algebra with central charge $c_{p,1}$. This is important
since the larger vertex operator algebra $M(1)$ is not completely
reducible.

 We prove, as the  main result, that
every irreducible $\overline{M(1)}$--module can be obtained as an
irreducible  subquotient of a certain $M(1)$--module $M(1,\l)$. In
order to prove this result, we determine explicitly  the Zhu's
algebra $A(\overline{M(1)})$. It is isomorphic to the commutative,
associative algebra ${\C}[x,y] / \la P(x,y) \ra$ where $\la P(x,y)
\ra $ is the  ideal in ${\C}[x,y]$ generated by  polynomial
\bea P(x,y) = y ^{2} - \frac{( 4 p) ^{2 p -1} }{  (2 p -1)! ^{2}}
\ ( x + \frac{(p-1) ^{2}}{4 p}) \prod_{i= 0} ^{p-2} \left( x +
\frac{i}{4 p} ( 2 p - 2 -i) \right) ^{2}. \nonumber \eea
This implies that the irreducible ${\N}$--graded
$\overline{M(1)}$--modules are parameterized by the solutions of
the equation $P(x,y) = 0$.
 The determination of the polynomial $P(x,y)$ is the central problem of our
construction. When $p=2$ (cf. \cite{W3}), this polynomial can be
constructed from a level-six singular vector in a generalized
Verma module over ${\W}(2,3)$--algebra. In general case, the
complicated structure of ${\W}(2, 2 p-1)$-algebras (cf.
\cite{EFH}, \cite{H}) makes also the calculation of  singular
vectors   extremely difficulty.

Since we are primary concentrated to the problem of classification
of irreducible representations, we only want to understand the
structure of the Zhu's algebra $A(\overline{M(1)})$. Fortunately,
this structure can be described without explicit knowledge of
relations among generators of  $ \overline{M(1)} $.

 In order to find
relations in  $A( \overline{M(1)} )$, we use a characterization of
the subalgebra $W_0$   as the kernel of the operator $Q \vert
M(1)$.   Using this, we show that the elements $H_{i} H$ for $i
\ge -2 p$ belong to the subalgebra $W_0$. It turns out that these
facts completely determine all relations in the Zhu's algebra
$A(\overline{M(1)})$ (cf. Section \ref{zhu-class}).

It is interesting, that the irreducible representations of the
Zhu's algebra $A(\overline{M(1)})$, and therefore of the vertex
operator algebra $\overline{M(1)}$ are parameterized by points of
a rational curve. Similar structures   were found for some other
irrational vertex operator (super)algebras (cf. \cite{A2}).

\section{Lattice and free boson vertex algebras }
\label{stwist}

We make the assumption that the reader is familiar with the
axiomatic theory of vertex (operator) algebras and  their
representations  (cf. \cite{DL}, \cite{FHL}, \cite{FZ},
\cite{FLM}, \cite{Z}).

In this section, we shall recall some properties of the lattice
and free boson vertex algebras.   The details can be found in
\cite{DL}, \cite{DMN}, \cite{K}, \cite{FLM}, \cite{MN}, \cite{S}.
Using language of generalized vertex operator algebras (cf.
\cite{DL}, \cite{GL}, \cite{S}), we construct screening operators
$ \widetilde{Q}$ and $Q$ acting on the vertex algebra $M(1)$.
Moreover, we choose the Virasoro element $\omega \in M(1)$ such
that $M(1)$ becomes a vertex operator algebra of rank $c_{p,1}$.
At the end of this section we shall present a result describing
the structure of $M(1)$ as a module for the Virasoro algebra with
central charge $c_{p,1}$.

\vskip 5mm

 Let $p \in {\N}$, $p \ge 2$. Let $\widetilde{L}={\Z}\beta$  be a
rational lattice of rank one with nondegenerate
 bilinear
 form $\la \cdot, \cdot \ra$ given by
 $$ \la \beta , \beta \ra = \frac{2}{p}.$$
 Let
 ${\h}={\C}\otimes_{\Z} \widetilde{L}$.
  Extend the form $ \la \cdot,
\cdot \ra $ on $\widetilde{L}$ to ${\h}$.
 Let $\hat{{\h}}={\C}[t,t^{-1}]\otimes {\h} \oplus {\C}c$ be the affinization of
${\h}.$
Set
$
\hat{{\h}}^{+}=t{\C}[t]\otimes
{\h};\;\;\hat{{\h}}^{-}=t^{-1}{\C}[t^{-1}]\otimes {\h}.
$
Then $\hat{{\h}}^{+}$ and $\hat{{\h}}^{-}$ are abelian subalgebras
of $\hat{{\h}}$. Let $U(\hat{{\h}}^{-})=S(\hat{{\h}}^{-})$ be the
universal enveloping algebra of $\hat{{\h}}^{-}$. Let ${\l} \in
{\h}$. Consider the induced $\hat{{\h}}$-module
\begin{eqnarray*}
M(1,{\l})=U(\hat{{\h}})\otimes _{U({\C}[t]\otimes {\h}\oplus
{\C}c)}{\C}_{\l}\simeq
S(\hat{{\h}}^{-})\;\;\mbox{(linearly)},\end{eqnarray*} where
$t{\C}[t]\otimes {\h}$ acts trivially on ${\C}$,
${\h}$ acting as $\la h, {\l} \ra$ for $h \in {\h}$
and $c$ acts on ${\C}$ as multiplication by 1. We shall write
$M(1)$ for $M(1,0)$.
 For $h\in {\h}$ and $n \in {\Z}$ write $h(n) =  t^{n} \otimes h$. Set
$
h(z)=\sum _{n\in {\Z}}h(n)z^{-n-1}.
$
Then $M(1)$ is a vertex   algebra which is generated by the fields
$h(z)$, $h \in {\h}$, and $M(1,{\l})$, for $\l \in {\h}$, are
irreducible modules for $M(1)$.
We shall choose the following Virasoro element in $M(1)$:
$$ \omega = \frac{p }{4 } {\beta} (-1) ^{2} + \frac{p-1}{2 } \beta
(-2).$$
The subalgebra of $M(1)$ generated by $\omega$ is isomorphic to
the simple vertex operator algebra $W_0=L(c_{p,1},0)$ where
$c_{p,1}=1 - 6 \frac{(p-1) ^{2}}{p}$. Let $$Y(\omega,z) = \sum_{n
\in {\Z} } L(n) z ^{-n -2} . $$
Thus $M(1)$ is a module for the Virasoro algebra (which we shall
denote by $Vir$)  with central charge $c_{p,1}$. In other words,
$M(1)$ becomes a Feigin-Fuchs module for the Virasoro algebra (cf.
\cite{FF}).

It is clear that $L(0)$ defines a $\N$--graduation on $M(1)=\oplus
_{ m \in {\N} } M(1)_m $. We shall write $\mbox{wt}(a) = m $ if $v
\in M(1)_m$. Thus $M(1)$ becomes a vertex operator algebra of rank
$c_{p,1}$ with the Virasoro element $\omega$.

Standard way for constructing vertex subalgebras of $M(1)$ is
given by the concept of screening operators acting on $M(1)$ (cf.
\cite{FFr}, \cite{FB}). We shall construct some particular
screening operators using language of generalized vertex operator
algebras (cf. \cite{DL}).

As in  \cite{DL}, \cite{S} (see also \cite{FLM}, \cite{K}), we
have the generalized vertex algebra
$$ V_{ \widetilde{L} } = M(1) \otimes {\C}[\widetilde{L}], $$
where ${\C}[\widetilde{L}]$ is a group algebra of $\widetilde{L}$
with a generator $e ^{\beta}$.
For $v \in  V_{ \widetilde{L} }$ let $Y(v,z) = \sum_{ s \in
\frac{1}{p} {\Z} } v_s z^{-s -1}$ be the corresponding vertex
operator (for precise formulae see \cite{DL}).

Clearly, $M(1)$ is a vertex subalgebra of $V_{ \widetilde{L} }$.

The Virasoro element $\omega \in M(1) \subset V_{ \widetilde{L} }$
is also a Virasoro element in $V_{\widetilde{L} }$ implying that
$V_{ \widetilde{L} }$ has a structure of a generalized vertex
operator algebra of rank $c_{p,1}$.

We have the  following decomposition:
$$ V_{\widetilde{L} } = \bigoplus_{ m \in {\Z} } M(1) \otimes e^{m
\beta}. $$

\begin{remark} {\rm If $p=2$ then $V_{  \widetilde{L} }$ is a
vertex operator superalgebra which can be constructed using
Clifford algebras (cf. \cite{K}, \cite{FB}). }
\end{remark}
\vskip 3mm

Define $\alpha = p \ \beta$. Then $\la \a , \a \ra = 2 p$,
implying that $L= {\Z} \alpha \subset \widetilde{L}$ is an even
lattice. Therefore the subalgebra $V_L \subset V_{ \widetilde{L}
}$ has a structure of a vertex operator algebra with the Virasoro
element $\omega$. In particular, for
 $v, w \in V_L$, we have  $Y(v,z) w= \sum_{n \in {\Z} } v_n w  \ z^{-n-1}$.

 Clearly,
$$ M(1) \subset V_L \subset V_{ \widetilde{L} }. $$
 Define the Schur polynomials $S_{r}(x_{1},x_{2},\cdots)$
 in variables $x_{1},x_{2},\cdots$ by the following equation:
\begin{eqnarray}\label{eschurd}
\exp \left(\sum_{n= 1}^{\infty}\frac{x_{n}}{n}y^{n}\right)
=\sum_{r=0}^{\infty}S_{r}(x_1,x_2,\cdots)y^{r}.
\end{eqnarray}

For any monomial $x_{1}^{n_{1}}x_{2}^{n_{2}}\cdots x_{r}^{n_{r}}$
we have an element $$h(-1)^{n_{1}}h(-2)^{n_{2}}\cdots
h(-r)^{n_{r}}{\1} $$ in $M(1)$   for $h\in{\h}.$
 Then for any polynomial $f(x_{1},x_{2}, \cdots)$,  $f(h(-1),
h(-2),\cdots){\1}$ is a well-defined element in $M(1)$ . In
particular, $S_{r}(h(-1),h(-2),\cdots){\1} \in M(1)$  for $r \in
{\N}$. Set $S_r (h)$ for $S_{r}(h(-1),h(-2),\cdots){\1}$.

 We shall now list some relations in the generalized vertex
operator algebra $V_{ \widetilde{L} }$.

Let $\gamma, \delta \in \widetilde{L}$. Instead of recalling  the
exact Jacobi identity in $V_{ \widetilde{L} }$,   we shall only
say that in the case  $$ \la \gamma, \delta \ra \in 2 {\Z}, $$ the
Jacobi identity  gives the following formulas (cf. \cite{DL},
\cite{GL}, \cite{S}):
\bea &&
 Y( e^{\gamma}, z) \ e^{\delta} = \sum_{ n \in {\Z} } e
^{\gamma} _n e ^{\delta} z ^{-n-1}, \\
&& \label{comut} [ e ^{\gamma} _n , e ^{\delta} _m] = \sum_{i=0}
^{\infty} { n \choose i } ( e ^{\gamma}_i e^{\delta}) _{ n + m -i}
. \ \ (m, n \in {\Z}) \eea

 The following relations in the generalized vertex
operator algebra $V_{ \widetilde{L} } $ are of great importance:
\begin{eqnarray}\label{eab1}
e ^{\gamma} _{i} e ^{\delta}=0\;\;\;\mbox{ for }i\ge
-\<\gamma,\delta\>.
\end{eqnarray}
Especially, if $\<\gamma,\delta\>\ge 0$, we have $e ^{\gamma} _{i}
e ^{\delta} =0$ for  $i\in {\N}$, and if
$\<\gamma,\delta\>=-n<0$,
 we get
 \bea\label{eab} && e ^{\gamma}
_{i-1} e ^{\delta} =S_{n-i}(\alpha) e ^{\gamma + \delta}
\;\;\;\mbox{ for }i\in \{ 0, \dots, n\}. \eea

 From the Jacobi identity in the (generalized) vertex operator
 algebras $V_L$ and $V_{ \widetilde{L}}$ follows:
\bea
&& [L(n), e^{\alpha} _m ] = -m   e ^{\alpha} _{m + n} ,
\label{kom}
\\
&&  [L(n), e ^{-\a} _m ] = ( 2 p  (n +1) - m) e ^{-\alpha} _ { n +
m}, \label{kom2}
\\
&& L(n) e ^{\alpha} = \delta_{n,0} e ^{\alpha} \ \ ( n \ge 0),
\label{kom3a}
\\ && L(n)  e ^{-\alpha} = \delta_{n,0} ( 2 p -1)
e ^{-\alpha} \ \ ( n \ge 0), \label{kom3} \\
&& L(n) e ^{-\beta} = \delta_{n,0} e ^{-\beta} \ (n \ge 0)
\label{kom-t-p} \\
&& [ L(n), e ^{-\beta} _r ] = - r e^{-\beta} _{r + n}, \ (r \in
\frac{1}{p} {\Z}). \label{kom-tilde}
 \eea

Define
\bea && Q =e ^{\a} _0 = \mbox{Res}_z   \ Y( e ^{\a},z) , \nonumber
\\ && \widetilde{Q} = e ^{-\beta}_0 = \mbox{Res}_z Y( e ^{-\beta},
z) \  . \nonumber \eea

From (\ref{kom}) and (\ref{kom-tilde}) we see that the operators
$Q$ and $ \widetilde{Q}$ commute with the Virasoro
 operators $L(n)$. We are interested in the action of these
 operators on $M(1)$. In fact, $Q$ and $ \widetilde{Q}$ are the screening
 operators, and therefore  $\mbox{Ker}_{M(1)} Q$ and
$\mbox{Ker}_{M(1)} \widetilde{Q}$ are vertex subalgebras of $M(1)$
(for details see Section 14 in \cite{FB} and reference therein).

Some properties of the screening operators $Q$ and $
\widetilde{Q}$ are given by the following lemma.

\begin{lemma} \label{rel-tilde} For  $ p >1$ we have:
\item[(i)] $[ Q , \widetilde{Q}] =0$.
\item[(ii)] $ \widetilde{Q} e ^{ n \alpha}  \ne 0$, $n \in
{\Zp}$.
\item[(iii)]$ \widetilde{Q} e ^{ -n \alpha}  =0$, $n \in
{\N}$.
\end{lemma} {\em Proof.} First we note that $ \la \a , \beta \ra = -2 \in
{\Z} $. Then the commutator formulae (\ref{comut}) gives that
$$[ Q, \widetilde{Q}] = [ e ^{\a} _0, e ^{-\beta} _0]= \left(
e^{\a} _0 e^{-\beta} \right)_{0}= \left( \a (-1) e^{\a - \beta}
\right)_0 = \frac{p}{p-1} \left(L(-1) e ^{\a - \beta} \right)_0.
$$
Since $(L(-1)u )_0 = 0$ in every generalized vertex operator
algebra, we conclude that $[ Q, \widetilde{Q}]=0$. This proves
(i).

Relation (ii)  follows from (\ref{eab}), and relation (iii) from
(\ref{eab1}). \qed

 We shall now investigate
  the action of the operator $Q$.
Since operators $Q ^{j}$, $j \in {\Zp}$, commute with the action
of the Virasoro algebra, they are  actually intertwiners between
Feigin-Fuchs modules inside the vertex operator algebra $V_L$.

Recall that a vector in $V_L$ is called primary if it is a
singular vector for the action of the Virasoro algebra.

Since $e ^{-n \a}$ is a primary vector in $V_L$
 for every $n \in {\N}$, we have that $Q ^{j} e ^{-n \a}$ is
 either zero or a primary vector.
Fortunately, the question of non-triviality of intertwiners
between Feigin-Fuchs modules is well studied in the literature.
So using arguments from \cite{F}, \cite{Ka1} together with the
methods developed in \cite{TK} and  \cite{Fel} one can see that
the following lemma holds.
\begin{lemma} \label{nontr-Q}
$ Q ^{j} e^{- n \alpha} \ne 0 \ \ \mbox{if and only if} \ \ j \in
\{ 0, \dots, 2 n \}. $
\end{lemma}

Next,  we shall present the theorem describing a structure of the
vertex operator algebra $M(1)$ as a module for the Virasoro vertex
operator algebra $L(c_{p,1},0)$. Again, the theorem can be proved
using a standard analysis in the theory of Feigin-Fuchs modules
(see \cite{FF}, \cite{Fel}, \cite{Ka1}, \cite{F}).

\begin{theorem} \label{str-ff}
 \item[(i)] The vertex operator algebra $M(1)$, as a module for the vertex
operator algebra $L(c_{p,1},0)$, is generated by the family of
singular and cosingular vectors $ \widetilde{Sing} \bigcup
\widetilde{CSing}$, where
$$  \widetilde{Sing} =  \{ u_{n} \ \vert \ n \in {\N}\}; \
\
  \widetilde{CSing} =  \{ w_{n} \ \vert \ n \in {\Zp} \}. $$
These vectors satisfy the following relations:
\bea
&& u_{n} = Q ^{n} e ^{-n \alpha} , \ \ Q ^{n} w_{n} = e ^{ n \a},
\nonumber \\
&& U(Vir) u_{n} \cong L(c_{p,1}, n ^{2} p + n p -n). \nonumber
\eea
\item[(ii)] The submodule generated by vectors $u_n, n \in {\N}$ is
isomorphic to $$ [\mbox{Sing}] \cong \bigoplus_{n =0 } ^{\infty}
L(c_{p,1}, n^{2} p + n p -n).$$

\item[(iii)]The quotient module  is isomorphic to
$$M(1)/  [\mbox{Sing}] \cong \bigoplus_{n =1 } ^{\infty}
L(c_{p,1}, n^{2} p - n p +n). $$
\item[(iv)] $Q u_0= Q {\1} =0$, and $Q u_n \ne 0$, $ Q w_n \ne 0$ for every $n
\ge 1$.
\end{theorem}

Theorem \ref{str-ff} immediately gives the following result.

\begin{proposition} \label{pomoc2}
 We have
$$L(c_{p,1},0) \cong  W_0 = \mbox{Ker}_{M(1)} Q  \ . $$
\end{proposition}

\section{The vertex operator algebra $ \overline{M(1)}$}  \label{def-voa-u}

  Recall that the Virasoro vertex operator
algebra $L(c_{p,1},0)$ is  the kernel of the screening operator
$Q$. But we have already seen that there is another screening
operator $ \widetilde{Q}$ acting on $M(1)$.
Define the following vertex algebra
$$\overline{M(1)} = \mbox{Ker}_{M(1) } \widetilde{Q}. $$

Since $ \widetilde{Q}$ commutes with the action of the Virasoro
algebra,   we have that
$$L(c_{p,1},0)\cong W_0 \subset  \overline{M(1)}. $$
This implies that $ \overline{M(1)} $ is  a vertex operator
subalgebra of $M(1)$ in the sense of \cite{FHL} (i.e.,  $
\overline{M(1)} $ has the same Virasoro element as $M(1)$).

The following theorem will describe the structure of the vertex
operator algebra $\overline{M(1)}$ as a $L(c_{p,1},0)$--module.

\begin{theorem} \label{generatori}
The vertex operator algebra $\overline{M(1)}$ is isomorphic to $
[\mbox{Sing}]$ as a $L(c_{p,1},0)$--module, i.e.,
$$ \overline{M(1)} \cong \bigoplus_{n =0} ^{\infty} L(c_{p,1},
n^{2} p + n p - n).$$
\end{theorem}
{\em Proof.} By Theorem \ref{str-ff} we know that the
$L(c_{p,1},0)$--submodule generated by the set $\widetilde{Sing} $
is completely reducible. So to prove the assertion, it suffices to
show that the operator $ \widetilde{Q}$ annihilates   vector $v
\in \widetilde{Sing} \cup \widetilde{CSing} $  if and only if $v
\in \widetilde{Sing} $.
Let $v \in  \widetilde{Sing} $, then $v = Q^{n} e ^{ -n \a}$ for
certain $n \in {\N}$. Since by Lemma \ref{rel-tilde}
$\widetilde{Q} e ^{ -n\a}=0$, we have that
$$ \widetilde{Q} v =\widetilde{Q} Q^{n} e ^{ -n \a} = Q^{n}
\widetilde{Q} e ^{ -n\a} = 0. $$

Let now $v \in \widetilde{CSing}$. Then there is $n \in {\Zp}$
such that $Q ^{n} v = e ^{ n \a}$. Assume that $ \widetilde{Q} v =
0$. Then we have that
$$0=  Q ^{n} \widetilde{Q} v = \widetilde{Q} Q ^{n} v =
\widetilde{Q} e ^{ n \a}, $$
contradicting Lemma \ref{rel-tilde} (ii). This proves the theorem.
 \qed

\begin{remark}
{\rm It is very interesting that although $M(1)$ is not completely
reducible $L(c_{p,1},0)$--module, its subalgebra $\overline{M(1)}$
is completely reducible. }
\end{remark}

Next  we shall prove that the vertex operator algebra $
\overline{M(1)} $ is generated by only two generators.

Motivated by formulae (18) in \cite{Ka1}, we define the following
three (non-zero) elements in $V_L$:
 $$ F =e ^{-\a}, \ \ H =   Q F, \ \ E =  Q ^{2}  F. $$
 From relations (\ref{kom}) - (\ref{kom3}) we see that
 $$ L(n) E = \delta_{n,0} (2 p  -1 ) E, \ L(n) F = \delta_{n,0}(2 p -1) F,
 \ L(n) H = \delta_{n,0}(2p-1) H \ \ ( n \ge 0 )$$
 i.e. $E$, $F$ and $H$ are primary vectors in $V_L$. In fact, $H$
 is a primary vector in $M(1)$.

 \begin{lemma} \label{pomoc1} In the vertex operator algebra $V_L$   the
 following relations hold:
\item[(i)] $ Q ^{3} F = 0$.
 \item[(ii)]  $E_{i} E = F_{i} F = 0$, for every $i \ge - 2 p$.
 \item[(iii)]  $ Q  ( H_{i} H )= 0$, for every $ i \ge - 2p$.
\item[(iv)] $H = S_{ 2 p -1} ( \alpha) $.
\end{lemma}
{\em Proof.}
Relation (i) is a special case of  Lemma \ref{nontr-Q}. Let now $i
\in {\Z}$, $i \ge - 2p$. From (\ref{eab1}) we have that $F_{i} F =
e ^ {-\alpha}_{i}  e ^{-\alpha} =0$.

Next we observe that $Q$ acts as a derivation on $V_L$, that is
\bea \label{q-der} Q ( a_n b) = (Q a)_n b + a_n (Q b) \ \ \
\mbox{for every} \ a, b \in V_L, \ n \in {\Z}.\eea

Then using (i)  and (\ref{q-der}) we see that $E_{i} E$ is
proportional to $Q^{4} (F_{i} F)$, which implies that $E_{i} E
=0$. This proves (ii).

Relation (iii) follows form (ii) and the fact that  $Q  ( H_{i} H
)$ is proportional to $ Q ^{3} ( F_{i} F)$.

Relation (iv) is a direct consequence of    (\ref{eab}).
 \qed \vskip 5mm

\begin{theorem}
The vertex operator algebra  $ \overline{M(1)} $ is  generated by
$\omega$ and $H=S_{2 p -1} (\alpha)$.
\end{theorem}
{\em Proof.}
Let $U$ be the vertex subalgebra of $ \overline{M(1)}$ generated
by $ \omega$ and $H$. We need to prove that $U = \overline{M(1)}$.
Let $W_n$ by the (irreducible) $Vir$--submodule of $
\overline{M(1)}$ generated by vector
$u_n$. Then $W_n \cong L(c_{p,1}, n ^{2} p + n p - n)$. %
Using Lemma \ref{nontr-Q} we see that
$$ \mbox{Ker}_{ \overline{M(1)} }\  Q ^{n} \cong \bigoplus_{i = 0
} ^{n -1 } W_i .$$
It suffices to prove that $u_n \in U$ for every $n \in {\N}$. We
shall prove this claim by induction. By definition we have that
$u_0, u_1(= H) \in U$. Assume that we have $k \in {\N}$ such that
$u_{n} \in U$ for $n \le k$. In other words, the inductive
assumption is
$ \oplus_{i = 0} ^{k} W_i \in U.$

We shall now prove that $u_{k+1} \in U$. Set $j = -2 k p -1$. By
Lemma \ref{nontr-Q} we have
$$ Q ^{2 k + 2} e ^{- (k+ 1) \a} =  Q ^{2 k + 2} \left(e ^{-\a}
_{j} e ^{- k \a} \right)  \ne 0. $$
Next we notice that
$$ Q ^{k + 1}  ( H _j u_k ) = Q ^{k + 1} \left( Q e^{-\a}\right)
_j \left( Q^{k} e^{-k\a} \right)= \frac{1}{2 k + 1}  Q ^{2 k + 2}
\left(e ^{-\a} _{j} e ^{- k \a} \right), $$
which implies that
 $$Q ^{k + 1}  (H _j u_k) \ne 0.$$
  So we have found vector $H_j u_k
 \in U$ such that
 $$ \mbox{wt} ( H_j u_k) = (2 p -1) + (k ^{2} p + k p - k) -j -1 =
 (k +1) ^{2} p  + (k+1) p - (k+1) = \mbox{wt} (u_{k+1}). $$
 This implies that
 $$ H_j u_k \in \bigoplus_{i=0} ^{k+1} W_i \ \ \mbox{and} \ \ \
  H_j u_k \notin \bigoplus_{i=0} ^{k} W_i \ .$$
 Since $Q ^{k+1 } \left(\oplus_{i=0} ^{k} W_i\right)= 0$ and
 $ \mbox{wt} ( H_j u_k) = \mbox{wt} (u_{k+1})$  we conclude that
 there is a   constant $C$, $C \ne 0$,  such that
 $$ H_j u _k = C u_{k + 1} + u', \ \ u' \in \bigoplus_{i=0} ^{k}
 W_i \subset U. $$
 Since $H_j u_k \in U$, we conclude that $u_{k+1} \in U$.

 Therefore, the claim is verified, and the proof of the theorem is complete. \qed

\begin{remark} {\rm
If $p=2$, then the elements $E, F, H$ span the triplet algebra
studied by M. Gaberdiel and H. Kausch (cf.  \cite{GK1},
\cite{GK2}, \cite{Ka2}). In this case  $ \overline{M(1)} $ is
isomorphic to ${\W}$--algebra ${\W}(2,3)$ with $c=-2$. Its
irreducible modules were classified by W. Wang in \cite{W2},
\cite{W3}.

In general, $E$, $F$ and $H$ span the  ${\W}(2, 2p-1, 2p-1, 2p
-1)$ algebra   with central charge $c_{p,1}$ (cf. \cite{Ka1},
\cite{F}). Our vertex operator algebra
  $ \overline{M(1)} $
 is isomorphic to the ${\W}$ algebra ${\W}(2, 2 p -1)$ investigated in a number of physical papers (cf. \cite{F},
\cite{EFH}, \cite{H}).

}
\end{remark}
\vskip 5mm
 The following lemma will imply that for $i \ge -2p$
vector $H_i H $ can be constructed using only the action of the
Virasoro operators $L(n)$ on  the vacuum vector ${\1}$.

\begin{lemma} \label{sing-z} We have:
 $$H_i H  \in W_0\cong L( c_{p,1}, 0) \ \ \mbox{ for every} \  i \ge - 2
 p. $$
In particular, $ H_{-1} H   \in W_0.$
\end{lemma}
{\em Proof.} The proof follows from Proposition \ref{pomoc2} and
Lemma \ref{pomoc1} (iii). \qed

\begin{remark}
{\rm In the case $p=2$, the fact that $H_{-1} H \in W_0$ can be
proved directly using a singular vector of level  $6 $ in a
generalized Verma module associated to ${\W}(2,3)$--algebra with
central charge $c=-2$ (cf. \cite{W3}, \cite{GK1}). This singular
vector implies that in  $ \overline{M(1)} $  the following
relation holds:
$$ H_{-1} H = \frac{1}{4}\left( \frac{19}{36} L(-3) ^{2} +
\frac{8}{9} L(-2) ^{3} + \frac{14}{9} L(-2) L(-4) - \frac{44}{9}
L(-6) \right) {\1} . $$

 Lemma \ref{sing-z}
indicates the existence of similar relation  of level $2 (2 p -1)$
in the general case, but we don't now the explicit formulae.
Instead of looking for such formulae, we use the realization of
the vertex operator algebra  $ \overline{M(1)} $ inside the
lattice vertex operator algebra $V_L$, and the description of the
Virasoro vertex operator algebra $L(c_{p,1},0)$ from Proposition
\ref{pomoc2}. }
\end{remark}

\section{ Spanning sets for $\overline{M(1)}$ and $A( \overline{M(1)})$ }

In this section shall  find a spanning set for  $ \overline{M(1)}
$ and for the Zhu's algebra $A( \overline{M(1)})$.

First we  recall   the definition of the Zhu's algebra for vertex
operator algebras.

Let $(V,Y, {\1}, \omega)$ be a vertex operator algebra. We shall
always assume that $V=\oplus_{ n \in {\N} } V_n$, where $V_n = \{
a \in V \ \vert \ L(0) a = n v \}$. For $a \in V_n$, we shall
write $\wt (a) = n$.

\begin{definition} {\rm We define the bilinear maps $* : V \otimes V
\rightarrow V$, $\circ : V \otimes V \rightarrow V$ as follows.
\bea
a*b &:= & \mbox{Res}_z  Y(a,z) \frac{ (1+z) ^{\wt (a) }}{z} b
= \sum_{i = 0} ^{\infty} { \wt (a) \choose i} a_{ i-1} b
, \nonumber \\
a\circ b &: = & \mbox{Res}_z  Y(a,z) \frac{ (1+z) ^{\wt (a) }}{z
^{2}} b
= \sum_{i = 0} ^{\infty} { \wt (a) \choose i} a_{ i-2} b
. \nonumber
 \eea
 Extend to $V \otimes V$ linearly, denote $O(V)\subset V$ the
 linear span of elements of the form $a \circ b$, and by $A(V)$
 the quotient space $V / O(V)$. }
\end{definition}
\vskip 5mm

Denote by $[a]$ the image of $a$ in $V$ under the projection of
$V$ into $A(V)$.
 We have:

 \begin{theorem} \cite{Z} \label{zhu-fund}
\item[(i)] The quotient space $(A(V),*)$ is an associative algebra
with unit element $[{\1}]$.

\item[(ii)] Let $M=\oplus_{ n \in {\N} } M(n)$ be a $\N$--graded
$V$--module. Then the top level $M(0)$ of $M$ is a $A(V)$--module
under the action $ [a]   \mapsto o(a) = a_{ \wt (a) -1}$ for
homogeneous $a \in V$.

\item[(iii)] Let $(U,\pi)$ be an irreducible $A(V)$--module. Then
there exists an irreducible ${\N}$--graded $V$--module
$L(U)=\oplus_{ n \in {\N} } L(U) (n)$ such that the top level
$L(U) (0)$ of $L(U)$ is isomorphic to ${U}$ as $A(V)$--module.

\item[(iv)] There is one-to-one correspondence between the
irreducible $A(V)$--modules and the irreducible ${\N}$--graded
$V$--modules.
 \end{theorem}

We shall need some information about the commutators $[H_n, H_m]$
for $m, n \in {\Z}$.

\begin{lemma} \label{komutatori}
For any $m, n \in {\Z}$, commutators $[H_n, H_m]$ are expressed as
(infinite) linear combination of
$$ L(p_1) \cdots L( p_s), \ \ p_1, \dots, p_s \in {\Z}, \ \ s \le
2 p -1. $$
In particular, for every vector $v \in \overline{M(1)}$ we have
$$ [H_n, H_m] v = f v, \ \ \mbox{for certain} \ f \in U(Vir) . $$
\end{lemma}
{\em Proof.}
 By the commutator formulae in vertex
(operator) algebras follows :
$$ [H_m, H_n] = \sum_{ i = 0} ^{\infty} { m \choose i} ( H_i H )_{
m + n -i} .$$
From Lemma \ref{sing-z} follows that $H_i H \in W_0$ for every
nonnegative integer $i$. Thus, $H_i H$ is an element of the
Virasoro vertex operator algebra $W_0 \cong L(c_{p,1},0)$. In
fact, $\mbox{wt} (H_i H) \le 2 (2 p -1)$, which implies that $H_i
H$ can be expressed as   (finite) linear combination
$$ L(-n_1 ) \cdots L( - n_s) {\1}, \ \ n_i \ge 2 ,  n_1 + \cdots +
n_s  \le 2 (2 p -1), \ s \le 2 p -1.$$
 Therefore, $[H_m, H_n]$ can be expressed as (infinite) linear
combination of
$$ L(p_1) \cdots L( p_s), \ \ p_1, \dots, p_s \in {\Z}, \ \ s \le
2 p -1. $$
This proves the first assertion. The second assertion follows from
the first assertion and from  the simple observation  that if $v
\in \overline{M(1)} $, and $m, n \in {\Z}$, then $[H_m , H_n ] v$
is well-defined. \qed

\begin{remark} {\rm
Using different arguments, a result  similar  to our Lemma
\ref{komutatori} was noticed  in the physical literature (cf.
\cite{EFH}, \cite{H}, \cite{KaW}).  }
\end{remark}

Lemma \ref{komutatori} shows that the structure of the vertex
operator algebra  $ \overline{M(1)} $ is similar
 to the structure of the vertex operator algebra $M(1)
^{+}$ studied in \cite{DN}.

So, using this lemma and  a completely analogous proofs to the
proofs of Proposition 3.4 and  Theorem 3.5 in \cite{DN}, one
obtains the following theorem.

\begin{theorem} \label{span-sets}
\item[(i)] The vertex operator algebra $\overline{M(1)}$ is spanned by
the following vectors
$$ L(-m_1) \cdots L(-m_s) H_{ -n_1} \cdots H_{ - n_t} {\1} , $$
where
$m_1 \ge m_2 \ge \cdots \ge m_s \ge 2$ and $ n_1 \ge n_2 \ge
\cdots \ge n_t \ge 1$.
\item[(ii)] The
Zhu's algebra $A( \overline{M(1)} )$ is spanned by the set
$$ \{ [ \omega ] ^{* s} * [H] ^{* t} \ \vert \ s, t \ge 0 \}. $$
In particular, the Zhu's algebra $A( \overline{M(1)})$ is
isomorphic to a certain quotient of the polynomial algebra ${\C}
[x, y]$, where $x$ and $y$ correspond $[\omega]$ and $[H]$.
\end{theorem}

 The fact that the Zhu's algebra $A(\overline{M(1)})$ is
commutative, enable us to study irreducible highest weight
representations of the vertex operator algebra $\overline{M(1)}$.
For given $(r,s) \in {\C} ^{2}$, let ${\C}_{r,s}$ be the one
dimensional module of $A( \overline{M(1)})$, with $[\omega]$
acting as the scalar $r$ and $[H]$ as the scalar $s$.  Therefore
every irreducible $A( \overline{M(1)})$--module is
one-dimensional, and it is isomorphic to a module ${\C}_{r,s}$ for
certain $(r,s) \in {\C} ^{2}$. Then Theorem \ref{zhu-fund} implies
that every irreducible ${\N}$--graded $\overline{M(1)}$--module is
isomorphic to a certain module $L({\C}_{r,s})$. So irreducible
representations of $\overline{M(1)}$ are parameterized by certain
subset of ${\C} ^{2}$. In the  following sections we will prove
that this subset is a rational curve.

\section{Representations of $\overline{M(1)}$ }
\label{representations}

In this section we identify a family of irreducible
$\overline{M(1)}$--modules. These modules are parameterized by
points $(r,s) \in {\C} ^{2}$ satisfying one algebraic equation.

 By construction, the vertex operator
algebra $\overline{M(1)}$ is a subalgebra  $M(1)$. We now that for
every $\lambda \in {\h} ( \cong {\C})$, $M(1,\lambda)$ is an
irreducible $M(1)$--module with the highest weight vector
$v_{\lambda}$. Thus $M(1,\lambda)$ is a $\overline{M(1)}$--module.
Denote by $\widetilde{V_{\lambda} }$ the
$\overline{M(1)}$--submodule generated by vector $v_{\l}$.

Set $ H(n) =   H_{ n + 2 p - 2}$, and
 $H(z) =  \sum_{ n \in {\Z} } H(n) z ^{-n - 2 p + 1}$.

First we recall the following result proved by the author in
\cite{A} for the purpose of studying ${\W}_{ 1 +
\infty}$--algebra.

\begin{proposition} \label{shur} [\cite{A}, Proposition 3.1]
Let $h \in {\h}$, and $r\in {\N}$.  Let $u =
S_{r}(h(-1),h(-2),\cdots){\1} $. Set $ Y( u,z) = \sum _{n \in {\Z}
} u_n z^{-n-1}$. Then we have
\bea (1) && u_n v_{\l} = 0 \quad \mbox{for} \ n>r-1, \nonumber \\
(2)&& u_{r-1} v_{\l} = { \la \l, h \ra \choose r} v_{\l}.
\nonumber \eea
\end{proposition}

Now Proposition \ref{shur} directly implies the following result.

\begin{proposition} \label{konstrukcija}
For every $\lambda \in {\h}$, $\widetilde{ {V}_{\lambda} }$ is a
${\N}$-graded $\overline{M(1)}$--module. The top level
$\widetilde{ {V}_{\lambda} }(0)$ is one-dimensional and generated
by $v_{\lambda}$. Let $ t = \la \lambda, \alpha \ra$. For every $n
\in {\N}$, we have
\bea  L(n) v_{ \lambda} &=& \delta_{n,0} \frac{1}{4 p} t ( t - 2
(p -1) ) v_{\lambda}, \label{tez-1} \\
H(n) v_{ \lambda} & = &  \delta_{n,0} { t \choose 2 p - 1} v_{
\lambda} \label{tez-2}. \eea
\end{proposition}

 By slightly abusing language, we can say that $\widetilde{V_{\l}}$ is a highest weight
$\overline{M(1)}$--module with respect to the Cartan subalgebra
$(L(0), H(0))$, and the highest weight is $(u(t) , v(t) )$ where
$t= \la \l, \alpha \ra \in {\C}$ and
\bea u(t) = \frac{1}{4 p} t ( t - 2 (p -1) ), \ \ v(t)=  { t
\choose 2 p - 1} \label{tez-3}.\eea
It is important to notice  that for every $t \in  {\C}$
$$ u(t) = u( 2(p-1) -t), \ \ v(t) = - v( 2( p-1) -t). $$
Moreover, the mapping $t \mapsto ( u(t), v(t) )$ is a injection.

Now we notice  that the top level of $\overline{M(1)}$--module
$\widetilde{V _{\lambda} }(0)={\C} v_{\l}$ has to be an
irreducible module for the Zhu's algebra $A(\overline{M(1)})$.
We have the following isomorphism of
$A(\overline{M(1)})$--modules:
$$ \widetilde{V _{\lambda} }(0) \cong {\C}_{r,s}, \ \
\mbox{where} \ \  r = u(t), \ s=v(t), \  t = \la \l , \a \ra. $$
Then one can show that the induced irreducible
$\overline{M(1)}$--module $L(\widetilde{V _{\lambda} }(0)) \cong
L({\C}_{r,s})$ is isomorphic to the irreducible quotient of $
\widetilde{V_{\l}}$.

Define $P(x,y) \in {\C}[x,y]$ by
%
\bea && P(x,y) = y ^{2} - C_p \  ( x + \frac{(p-1) ^{2}}{4 p})
\prod_{i= 0} ^{p-2} \left( x + \frac{i}{4 p} ( 2 p - 2 -i) \right)
^{2}, \label{def-pol} \eea
where $C_p =    ( 4 p) ^{2 p -1} /  (2 p -1)! ^{2}$.

In the following lemma we shall see that the highest weights of
  $ \overline{M(1)}$--modules $ \widetilde{V_{\l}}$ $(\l \in {\h})$
   coincide with the solutions
  of the equation $P(x,y) =0$. The
 proof is similar to the proof of Lemma 4.4 in \cite{W3}.

\begin{lemma} \label{parametrizacija}
Solutions of the equation
\bea && P(x,y) = 0 \label{tez-rel} \eea
are parameterized by
\bea (x, y) = (u(t), v(t)), \ \ t \in {\C}. \label{para} \eea
\end{lemma}
{\em Proof.} It is easy to  verify that for every $t \in {\C}$
$P(u(t), v(t)) =0$.

Let now $(x,y)$ be any solution of (\ref{tez-rel}). Then there is
$t_0 \in {\C}$ such that
$$ x = u( t_0) = u ( 2(p - 1) - t_0). $$
By substituting $x= u(t_0)$ in the equation (\ref{tez-rel}) we get
that
$y ^{2} = { t_0 \choose 2 p -1} ^{2}$, which implies that
$$ y =  { t_0 \choose 2 p -1} = v(t_0) \ \mbox{or} \ \ y = - { t_0
\choose 2 p -1} = v(2(p -1) -t_0).$$
So there is a unique $t \in {\C}$ such that (\ref{para}) holds.
\qed

\begin{theorem} \label{zhu-m-konst}
Assume that $(r,s) \in {\C} ^{2}$ so that $P(r,s) = 0$. Then
\item[(i)] ${\C}_{r,s}$ is an irreducible $A(
\overline{M(1)})$--module.
\item[(ii)] $L({\C}_{r,s})$ is an irreducible
$\overline{M(1)}$--module.
\end{theorem}
{\em Proof.}  Using Lemma \ref{parametrizacija} we see that there
is a unique $\l \in {\h}$ so that ${\C}_{r,s} \cong
\widetilde{V_{\l}}(0) $, where  $\widetilde{V_{\l}}(0) $ is  an
$A( \overline{M(1)})$--module constructed in Proposition
\ref{konstrukcija}. So ${\C}_{r,s}$ is an irreducible $A(
\overline{M(1)})$--module. Assertion (ii) follows from Theorem
\ref{zhu-fund}. \qed

\section{ The Zhu's algebra $A( \overline{M(1)}  )$ and the
classification of irreducible modules} \label{zhu-class}

In this section we shall prove our main result saying that the
modules constructed in Section \ref{representations} provide all
irreducible ${\N}$--graded $\overline{M(1)}$--modules. Our proof
will use the theory  of Zhu's algebras.

 We are now going to determine  the
Zhu's algebra for the vertex operator algebra $\overline{M(1)}$.
We have already proved that $A( \overline{M(1)} )$ is isomorphic
to a certain quotient of the polynomial algebra ${\C}[x,y]$, where
$x$ corresponds to $[\omega]$ and $y$ to $[H]$. Now we shall find
all relations in $A( \overline{M(1)} )$.

\begin{lemma} \label{rel-zhu}
  In the Zhu's algebra $A( \overline{M(1)})$, we have:
$$ P ( [\omega], [H] ) = 0. $$
\end{lemma}
{\em Proof.}
 Lemma \ref{sing-z} implies that for every $i > - 2 p$
$$ H_{i-1} H   \ = \  f_{i} ( L(-2), L(-3), \dots  ) {\1}$$
for certain polynomial $f_i  \in {\C} [x_1, x_2, \dots ]$.
This implies that in  $A( \overline{M(1)})$, we have
$$[H_{i-1} H] = g_i ([L(-2){\1}]) = g_i ([\omega])$$
for certain polynomial $g_i \in {\C}[x]$ such that $\deg (g_i) \le
2 p -1$.
The definition of the multiplication in $A( \overline{M(1)})$
gives that
$$ [H] * [H] = \sum_{ i = 0} ^{ 2 p -1 } { 2 p -1 \choose i } [ H_
{i-1} H]=  \sum_{ i = 0} ^{ 2 p -1 } { 2 p -1 \choose i } g_i(
[\omega]). $$
Let  $g(x)= \sum_{ i = 0} ^{ 2 p -1 } { 2 p -1 \choose i }
g_i(x)$.
So we have proved that there is  polynomial $g \in {\C}[x]$ so
that
\bea \label{rel-ev} [H] * [H] = [H] ^{2} = g ([\omega]),  \ \deg
(g) \le 2 p -1. \eea
Now we shall  determine the polynomial $g$ explicitly.
Recall that if $(r,s) \in {\C} ^{2}$ such that $P(r,s)=0$, then
${\C}_{r,s}$ is an irreducible $A( \overline{M(1)})$--module
(Theorem \ref{zhu-m-konst}). Let us now evaluate both sides of
(\ref{rel-ev}) on  $A( \overline{M(1)})$--modules ${\C}_{r,s}$. We
get that
$ s ^{2} = g(r) $ for every $(r,s) \in {\C} ^{2}$ such that
$P(r,s)=0$. This implies that for every $r \in {\C}$
 $$ g(r) =s ^{2}= s^{2} - P(r,s) =C_p \  ( r + \frac{(p-1) ^{2}}{4 p}) \prod_{i= 0}
^{p-2} \left( r + \frac{i}{4 p} ( 2 p - 2 -i) \right) ^{2}, $$
where $C_p =   ( 4 p) ^{2 p -1} / (2 p -1)! ^{2}$.
In this way we have proved that
 $$ g(x) =  y ^{2} - P(x,y)  =C_p \  ( x + \frac{(p-1) ^{2}}{4 p}) \prod_{i= 0}
^{p-2} \left( x + \frac{i}{4 p} ( 2 p - 2 -i) \right) ^{2} .
  $$
Using (\ref{rel-ev}) we get
$$ P( [\omega], [H] ) =  [H] ^{2} - g ([\omega]) = 0, $$
as desired.
 \qed
\vskip 5mm  Now we are in the position to find all relations in
the Zhu's algebra $A( \overline{M(1)} )$.
\begin{theorem} \label{zhu-alg}
The Zhu's algebra $A(\overline{M(1)})$ is isomorphic to the
commutative, associative algebra  ${\C}[ x, y] / \la P(x,y) \ra$.
\end{theorem}
{\em Proof.}
By Theorem \ref{span-sets} we have a surjective algebra
homomorphism
\bea
\Phi :  {\C}[x,y]&& \rightarrow A( \overline{M(1)} ) \nonumber \\
 x &&\mapsto [\omega] \nonumber \\
y &&\mapsto [H] . \nonumber \eea
It suffices to prove that $\mbox{Ker} \  \Phi = \la P(x,y) \ra$.

 Lemma \ref{rel-zhu} gives  that $ \la P(x,y) \ra \subseteq \mbox{Ker}
\  \Phi$.

Assume now that $K(x,y) \in \mbox{Ker} \ \Phi$. Note that $P(x,y)$
has degree $2$ in $y$. Using the division algorithm we  get
$$ K(x,y) = A(x,y) P(x,y) + R(x,y), $$
where $A(x,y), R(x, y) \in {\C} [x,y]$ so that $R(x,y)$ has degree
at most $1$ in $y$. So we can write
$ R(x,y) = B(x) y + C(x)$, where $B(x) , C(x) \in {\C}[x]$.  Since
$P(x,y), K(x,y) \in \mbox{Ker} \ \Phi$ we have that $R(x,y) \in
\mbox{Ker} \ \Phi$.  We  now evaluate polynomial $R(x,y)$ on $A(
\overline{M(1)})$--modules and obtain
$$ R( u(t), v(t) ) = 0 \ \mbox{for every} \ t \in {\C},$$
where polynomials $u(t), v(t)$ are defined by (\ref{tez-3}).
Therefore
\bea \label{rel-pom3}
B( u(t) ) v(t) = - C(t) \ \ \mbox{for every} \ t \in {\C}. \eea
Assume that $B(x)\ne 0$. Then polynomial $$B( u(t) ) v(t) = B(
\frac{1}{4p} ( t ( t- 2 p +2)  ) {t \choose 2 p -1}$$
has odd degree in $t$. On the other hand, if  polynomial $- C(
u(t) )$ is  nontrivial, it must have even degree in $t$. This is a
contradiction because of (\ref{rel-pom3}). So $ B(x)  = 0$. Using
(\ref{rel-pom3}) we also get that $C(u(t) )=0 $ for every $t \in
{\C}$, which implies   that $C(x)=0$. In this way we have proved
that
$$ R(x,y) = B(x) y + C(x) = 0 , $$
and therefore
$$ K(x,y) = P(x,y) A(x,y) \ \in \ \la P(x,y) \ra .$$
So $ \mbox{Ker} \  \Phi = \la P(x,y) \ra$, and the theorem holds.
\qed

\begin{theorem} \label{class-p}
The set
\bea \label{skup-svih} \{  L( {\C}_{r,s} )   \ \vert \  (r,s) \in
{\C} ^{2}, \ P(r,s) = 0 \} \eea provides all non-isomorphic
irreducible $\N$--graded modules for the vertex operator algebra
$\overline{M(1)}$.
\end{theorem}
{\em Proof.} Since the Zhu's algebra $A(\overline{M(1)})$ is
commutative,  Theorem \ref{zhu-alg} implies that the set
$$ \{ {\C}_{r,s} \ \ \vert \ \ (r,s) \in {\C} ^{2}, \ P(r,s)=0 \
\}
$$ provides all  irreducible modules for the Zhu's algebra $A(
\overline{M(1)})$. Then Theorem \ref{zhu-fund} implies that the
set (\ref{skup-svih}) provides all irreducible ${\N}$--graded
$\overline{M(1)}$--modules.  \qed

\begin{remark}{\rm  Theorem \ref{class-p} shows that the
 irreducible  $\overline{M(1)}$--modules  are parameterized by
the solutions of the equation $P(x,y) =0$. We have observed that
all solutions of this equation can be written in the form
$( u(t), v(t) )$, $(t \in {\C})$ which are exactly the highest
weights of $ \overline{M(1)}$--submodules of $M(1,\l)$ constructed
in Proposition \ref{konstrukcija}. This leads to the  conclusion
(as in \cite{W3}) that every irreducible $\overline{M(1)}$--module
can be identified starting from modules for the   vertex operator
algebra $M(1)$.
}
\end{remark}

\end{document}